\documentclass[11pt, reqno]{amsart}
\usepackage{amssymb}
\usepackage{amscd}
\usepackage{verbatim,ifthen}
\usepackage{color}
\usepackage{latexsym}
\usepackage{tikz}
\usepackage{tikz-cd}
\usepackage{mathrsfs}
\usepackage{wrapfig}
\usetikzlibrary{shapes}
\usepackage{color}
\usetikzlibrary{arrows.meta}
\usepackage{bbm}
\usetikzlibrary{matrix}
\usetikzlibrary{calc}
\usetikzlibrary{arrows,intersections}
\usepackage{pgfplots}
\usepackage{multicol}
\usepackage{array}
\newcolumntype{M}[1]{>{\centering\arraybackslash}m{#1}}

\usepackage[colorlinks, linkcolor=black, citecolor=magenta, linktocpage]{hyperref}
\title{Twisted K\"ahler--Einstein Metrics and Collapsing}
\author{Kyle Broder }
\addtocontents{toc}{\setcounter{tocdepth}{1}}

\let\oldtocsection=\tocsection
 
\let\oldtocsubsection=\tocsubsection
 
\let\oldtocsubsubsection=\tocsubsubsection
 
\renewcommand{\tocsection}[2]{\hspace{0em}\oldtocsection{#1}{#2}}
\renewcommand{\tocsubsection}[2]{\hspace{1em}\oldtocsubsection{#1}{#2}}
\renewcommand{\tocsubsubsection}[2]{\hspace{2em}\oldtocsubsubsection{#1}{#2}}

\def\XXint#1#2#3{{\setbox0=\hbox{$#1{#2#3}{\int}$ }
\vcenter{\hbox{$#2#3$ }}\kern-.6\wd0}}

\usepackage{eucal}
\usepackage{calc}  
\usepackage{enumitem} 
\usepackage{tensor}
\usepackage{graphicx,wrapfig,lipsum}
\usepackage{etoolbox}
\usepackage{marginnote}
\usepackage{lipsum}
\makeatletter
\patchcmd{\@mn@margintest}{\@tempswafalse}{\@tempswatrue}{}{}
\patchcmd{\@mn@margintest}{\@tempswafalse}{\@tempswatrue}{}{}
\reversemarginpar 
\makeatother
\usepackage{scrextend}

\makeatletter
\DeclareRobustCommand\widecheck[1]{{\mathpalette\@widecheck{#1}}}
\def\@widecheck#1#2{%
    \setbox\z@\hbox{\m@th$#1#2$}%
    \setbox\tw@\hbox{\m@th$#1%
       \widehat{%
          \vrule\@width\z@\@height\ht\z@
          \vrule\@height\z@\@width\wd\z@}$}%
    \dp\tw@-\ht\z@
    \@tempdima\ht\z@ \advance\@tempdima2\ht\tw@ \divide\@tempdima\thr@@
    \setbox\tw@\hbox{%
       \raise\@tempdima\hbox{\scalebox{1}[-1]{\lower\@tempdima\box
\tw@}}}%
    {\ooalign{\box\tw@ \cr \box\z@}}}
\makeatother

\begin{document}

\maketitle

\begin{abstract}
We extend the arguments of Tosatti--Zhang \cite{TosattiZhang} to reduce a well-known conjecture concerning the structure of the Gromov--Hausdorff limit in both the setting of degenerating Calabi--Yau manifolds and the K\"ahler--Ricci flow to a certain partial second-order estimate. 
\end{abstract}

\section{Introduction}
Let $M$ be a compact Riemannian manifold without boundary. An old problem asks whether $M$ admits a canonical metric. There are many variants of this canonical metric problem. If $M$ happens to be K\"ahler, then there has been great success in finding K\"ahler--Einstein metrics: that is, K\"ahler metrics with constant Ricci curvature. The existence of these metrics, however, imposes very restrictive conditions on the topology of such manifolds. In particular, K\"ahler--Einstein manifolds necessarily have a definite, or trivial, first Chern class. 

Hence, one is naturally inclined to ask whether a canonical metric can exist on a compact K\"ahler manifold which does not fit into one of these restrictive categories: $c_1(M) >0$, $c_1(M) = 0$, or $c_1(M) <0$. Such metrics were introduced in \cite{SongTian2007, SongTian2012} as long-time solutions of the K\"ahler--Ricci flow, and also arise in the study of collapsing degenerations of Calabi--Yau metrics \cite{TosattiAdiabaticLimits}. In both of these settings there is a holomorphic fibration\footnote{That is, a surjective holomorphic map with connected fibers.} $f : M \longrightarrow N$ and a canonical metric $\omega_{\text{can}}$ is known to exist on the base of this fibration.\footnote{See \cite[Definition 3.4]{SongTian2012} for a more precise definition of canonical metric.} Moreover, on the regular locus $N_0$ of this fibration\footnote{By which we mean $N_0 : = N \backslash D$, where $D$ is the discriminant locus of $f$, consisting of the singular locus of $N$ together with the critical values of $f$.}, the canonical metric satisfies the twisted K\"ahler--Einstein equation $$\text{Ric}(\omega_{\text{can}}) = \lambda \omega_{\text{can}} + \omega_{\text{WP}},$$ where $\lambda \in \{ -1, 0 \}$ (depending on the setting: $\lambda =0$ in the Calabi--Yau case, while $\lambda =-1$ for the K\"ahler--Ricci flow) and $\omega_{\text{WP}}$ is the Weil--Petersson metric pulled back from the moduli space of Calabi--Yau manifolds (c.f., \cite{FangLu}), measuring the change of complex structure in the fibers of $f$. 

Let $D \subset N$ denote the discriminant locus of $f$, i.e., the set of points in $N$ over which the fibers of $f$ are singular. We set $N_0 : = N \backslash D$. Let $\pi : \widetilde{N} \longrightarrow N$ be a log resolution of the pair $(N,D)$. In particular, $\pi$ is a surjective holomorphic map with connected fibers, $\widetilde{N}$ is smooth, and $E = \pi^{-1}(D)$ has simple normal crossings. It was discovered in \cite{TosattiZhang} that certain estimates on $\pi^{\ast} \omega_{\text{can}}$ imply global results concerning the structure of the Gromov--Hausdorff limit of degenerating Calabi--Yau metrics and long-time solutions of the K\"ahler--Ricci flow. We will consider both settings here. Before considering both cases in more detail, let us state the main theorem of the paper:

\subsection*{Main Theorem}\label{Main Theorem}
Let $(X, d_X)$ denote the metric completion of $(N_0, \omega_{\text{can}})$, and set $S_X = X - N_0$. Suppose that for any $\varepsilon > 0$, there exists a constant $C >0$ (possibly depending on $\varepsilon$) and $d \in \mathbb{N}$ such that on $\widetilde{N} \backslash E$, \begin{eqnarray}\label{MAINESTIMATE1}
\pi^{\ast} \omega_{\text{can}} & \leq & \frac{C}{| s_F |^{2\varepsilon}} \left( 1 - \sum_{i=1}^{\mu} \log | s_i |_{h_i} \right)^d \omega_{\text{cone}},
\end{eqnarray}
then, in the Calabi--Yau case, \begin{itemize}
\item[(i)] $(X, d_X)$ is a compact length metric space and $S_X$ has real Hausdorff codimension at least 2. 
\item[(ii)] $(M, \widetilde{\omega}_t)$ converges in the Gromov--Hausdorff topology to $(X, d_X)$ as $t \to 0$. 
\item[(iii)] $X$ is homeomorphic to $N$. 
\end{itemize}
Moreover, assuming a uniform bound on the Ricci curvature on compact subsets of $M \backslash S$, (i)--(iii) also hold for the K\"ahler--Ricci flow. \\

This weakens the estimate proposed in \cite{TosattiZhang} to prove (i)--(iii). Indeed, there it is proposed that the multiplicity of the divisorial pole (i.e., $2\varepsilon$) should be zero. Here we show that it suffices to have $\varepsilon >0$ small enough.

\subsection{Collapsing Degenerations of Calabi--Yau Metrics}
Let $(M^m, \omega_M)$ be a compact Calabi--Yau manifold of dimension $m$, which for simplicity, we will assume is projective, i.e., $M$ is a projective manifold with $K_M \simeq \mathcal{O}_M$. We suppose that there is a fibration $f : M \longrightarrow N$ over a normal projective variety $N$ of dimension $n < m$. As before, we let $D$ denote the discriminant locus of $f$. Setting $S = f^{-1}(D)$, it follows that $f : M \backslash S \longrightarrow N \backslash D$ is a proper holomorphic submersion. Let $\omega_N$ be a K\"ahler metric on $N$ (understood in the sense of analytic spaces for non-smooth $N$, see, e.g., \cite{Moishezon, EGZKE}). Then by \cite{SongTian2012, TosattiAdiabaticLimits,EGZKE,CGZ} it is known that there exists a continuous $\omega_N$--plurisubharmonic function $\varphi$ such that $\omega_{\text{can}} : = \omega_N + \sqrt{-1} \partial\overline{\partial} \varphi$ satisfies the twisted K\"ahler--Einstein equation $$\text{Ric}(\omega_{\text{can}}) \ = \ \omega_{\text{WP}} \ \geq \ 0,$$ on $N_0 = N \backslash D$. Note that $\omega_{\text{WP}} \equiv 0$ if the fibers of $f$ are all biholomorphic. 

We want to understand the nature of the singularities of $\omega_{\text{can}}$ near $D$. The first general result in this direction is due to \cite{Hein2012}, where very precise estimates are obtained when the base of the fibration is a Riemann surface. Later results were obtained in \cite{TosattiZhang} when $M$ is hyperk\"ahler. Very recently, these results were generalized and improved in \cite{GTZ2019}.

We recall some results from \cite{GTZ2019}. Let $N^{\text{reg}}$ denote the regular locus of $N$ and decompose $D = D^{(1)} \cup D^{(2)}$ into a codimension-one part $D^{(1)}$ and a codimension-two part $D^{(2)}$. Let $D^{\text{snc}} \subset D^{(1)} \cap N^{\text{reg}}$ be the simple normal crossings locus of $D$. Decompose $D = \bigcup_{i=1}^{\mu} D_i$ into irreducible components and let $s_i$ be the defining sections of the line bundles $\mathcal{O}(D_i)$ associated to the divisor $D_i$. Endow the line bundles $\mathcal{O}(D_i)$ with Hermitian metrics $h_i$. As detailed in $\S 3$, we can associate a conical K\"ahler metric $\omega_{\text{cone}}$ to the snc data $\{ D_i, 2\pi \alpha_i \}$, where $\alpha_i \in (0,1]$ are the cone angles of $\omega_{\text{cone}}$ along $D_i$.

\subsection*{Theorem 1.1.1}
(\cite[Theorem 1.1]{GTZ2019}). The canonical twisted K\"ahler--Einstein metric $\omega_{\text{can}}$ extends smoothly across $D^{(2)} \cap N^{\text{reg}}$ and there exists a conical K\"ahler metric $\omega_{\text{cone}}$ as above such that for any point $x \in D^{\text{snc}}$, there is an open set $U \subset N^{\text{reg}}$ and constants $d \in \mathbb{N}$, $C>0$ such that \begin{eqnarray*}
C^{-1} \left( 1 - \sum_{i=1}^{\mu} \log | s_i |_{h_i} \right)^{-d\max(n-2,0)} \omega_{\text{cone}} \ \leq \ \omega_{\text{can}} \ \leq \ C \left( 1 - \sum_{i=1}^{\mu} \log | s_i |_{h_i} \right)^d \omega_{\text{cone}}.
\end{eqnarray*}

\vspace*{0.2cm}

The logarithmic terms are negligible in comparison with the poles of $\omega_{\text{cone}}$. This is an important observation for later applications to collapsing. Let $\pi : \widetilde{N} \longrightarrow N$ be as before. Again, we can associate a conical K\"ahler metric $\omega_{\text{cone}}$ to the snc data $\{ E_i, 2\pi \alpha_i \}$, where $E = \bigcup_i E_i$ is a decomposition of $E$ into irreducible components. Let $F = \bigcup_j F_j \subset E$ is the union of $\pi$--exceptional components of $E$.

\subsection*{Theorem 1.1.2}\label{GTZTheorem}
(\cite[Theorem 5.2]{GTZ2019}). On $\widetilde{N} \backslash E$, there exists a uniform consant $C>0$ and $d \in \mathbb{N}$ such that \begin{eqnarray}\label{GTZMAIN}
\pi^{\ast} \omega_{\text{can}} \ \leq \ \frac{C}{| s_F |^{2A}} \left( 1- \sum_{i=1}^{\mu} \log |s_i |_{h_i} \right)^d \omega_{\text{cone}}.
\end{eqnarray}

\vspace*{0.2cm}

where $| s_F |^2 : = \prod_j | s_{F_j} |_{h_j}^2$. \\

We note that $A$ has to be taken very large in the proof of \cite[Theorem 5.2]{GTZ2019}. This prevents the arguments in \cite{TosattiZhang} from going through. 

\subsection*{Corollary 1.1.3}
In the special case that $N$ is smooth and $D^{(1)}$ is a simple normal crossings divisor, $$\pi^{\ast} \omega_{\text{can}} \ \leq \ C \left( 1 - \sum_{i=1}^{\mu} \log | s_i |_{h_i} \right)^d \omega_{\text{cone}}.$$

\vspace*{0.2cm} 

Indeed, under the assumption that $N$ is smooth and $D^{(1)}$ is snc, the resolution $\pi$ is simply the identity map and there is no divisorial terms in \eqref{GTZMAIN}. In light of these results and the main theorem here, we conjecture the following partial second-order estimate:

\subsection*{Conjecture 1.1.4}\label{MainTheorem}
For any $\varepsilon > 0$, there exists a constant $C>0$ (dependent on $\varepsilon$) and $d \in \mathbb{N}$ such that on $\widetilde{N} \backslash E$, \begin{eqnarray}\label{PARTIALSECOND}
\pi^{\ast} \omega_{\text{can}} \ \leq \ \frac{C}{| s_F |^{2\varepsilon}} \left( 1- \sum_{i=1}^{\mu} \log |s_i |_{h_i} \right)^d \omega_{\text{cone}}.
\end{eqnarray}

\vspace{0.2cm}

\subsection{Applications to Collapsing Calabi--Yau Manifolds}
Let $M^m$ be a Calabi--Yau $m$--fold which is the total space of a fibration $f : M^m \longrightarrow N^n$ over some normal projective variety $N$ of dimension $n<m$. Let $$\mathcal{K} \ : = \ \left \{ \alpha \in H^{1,1}(M, \mathbb{R}) : \exists \ \omega \ \text{K\"ahler on $M$ such that} \ [\omega] =\alpha \right \}$$ be the K\"ahler cone of $M$.\footnote{This is an open convex cone in the finite-dimensional vector space $H^{1,1}(M, \mathbb{R})$.} Let $\alpha_0$ be a fixed non-zero class on the boundary of $\mathcal{K}$, which we will take to be $\alpha_0 = f^{\ast}[\omega_N]$, for some K\"ahler metric $\omega_N$ on $N$. For $0 < t \leq 1$, we consider the path $\alpha_t = f^{\ast}[\omega_N] + t[\omega_M]$ which approaches $f^{\ast}[\omega_N]$ from the interior of the K\"ahler cone as $t \to 0$. By Yau's solution of the Calabi conjecture \cite{Yau1976}, there is a unique Ricci-flat K\"ahler metric $\widetilde{\omega}_t$ in each class $\alpha_t$ for $t >0$.

These metrics $$\widetilde{\omega}_t = f^{\ast} \omega_N + t \omega_M + \sqrt{-1} \partial \overline{\partial} \varphi_t, \hspace*{1cm} \sup_M \varphi_t =0,$$ satisfy the complex Monge--Amp\`ere equation $$\widetilde{\omega}_t^m = c_t t^{m-n} \omega_M^m,$$ where the constants $c_t$ are bounded away from $0$ and $\infty$ and converge as $t \to \infty$. The convergence of these metrics has been a well-studied problem. In \cite{TosattiAdiabaticLimits}, it was shown that $\widetilde{\omega}_t$ converges to the (pullback of the) twisted K\"ahler--Einstein metric $f^{\ast}\omega_{\text{can}}$ in the $\mathcal{C}_{\text{loc}}^{1,\gamma}(M \backslash S)$ topology of K\"ahler potentials for all $0 < \gamma <1$. This was improved to local uniform convergence in \cite{TosattiWeinkoveYang} and to $\mathcal{C}_{\text{loc}}^{\alpha}(M \backslash S)$ in \cite{HeinTosatti}. The convergence occurs in $\mathcal{C}_{\text{loc}}^{\infty}(M \backslash S)$ if the generic fibers are tori, or finite quotients of tori \cite{GTZ2013,HeinTosattiRemarks,TosattiZhangInfiniteTimeSingularities}, or if the generic fibers are all biholomorphic \cite{HeinTosatti}. 

One of the main conjectures that has been of interest in recent years is the following:

\subsection*{Conjecture 1.2.1}\label{MainConjecture}
Let $(X, d_X)$ be the metric completion of $(N_0, \omega_{\text{can}})$. Set $S_X = X - N_0$. Then \begin{itemize}
\item[(i)] $(X, d_X)$ is a compact length metric space and $S_X$ has real Hausdorff codimension at least 2. 
\item[(ii)] $(M, \widetilde{\omega}_t)$ converges in the Gromov--Hausdorff topology to $(X,d_X)$ as $t \to 0$.
\item[(iii)] $X$ is homeomorphic to $N$. 
\end{itemize}

\vspace*{0.2cm}

\nameref{MainConjecture} is motivated by analogous conjectures made in \cite{GrossWilson, KontsevichSoibelman, Todorov} for collapsed limits of Calabi--Yau manifolds near a large complex structure limit. The conjecture was first verified by \cite{GrossWilson} when $f : M \to N$ is an elliptic fibration of K3 surfaces with only $I_1$ singular fibers. This uses a precise gluing construction that is difficult to generalize in higher-dimensions due to the large number of possible types of singular fibers (see also \cite{ChenViaclovskyZhang, Li}). When $N$ is a Riemann surface \nameref{MainConjecture} was verified by \cite{GTZ2016}. 

When $K_M$ is nef and big, the Gromov--Hausdorff limit was identified in \cite{BingWang}. In \cite{SongTianZhang2019}, part (ii) of \nameref{MainConjecture} was proved in general, while part (iii) was established when $N$ has at worst orbifold singularities. Finally, following the program laid out in \cite{TosattiZhang}, the conjecture was recently verified by \cite{GTZ2019} assuming that $D^{(1)}$ is a simple normal crossings divisor. 

\subsection{The K\"ahler--Ricci flow}
Fix a compact K\"ahler manifold $M^m$ of dimension $m$. In \cite{SongTian2007, SongTian2012} an \textit{Analytic Minimal Model Program} was proposed to study the birational classification of compact K\"ahler manifolds (in particular, algebraic manifolds) using the K\"ahler--Ricci flow:

\begin{eqnarray}\label{KRF}
\frac{\partial \omega(t)}{\partial t} \ = \ - \text{Ric}(\omega(t)) - \omega(t), \hspace*{1cm} \omega(0) = \omega_0,
\end{eqnarray} where $\omega_0$ is some fixed K\"ahler metric on $M$.

In \cite[$\S 6.2$]{SongTian2017} a series of conjectures were made on the behavior of the K\"ahler--Ricci flow, together with their implications on birational classification. In short, it is predicted that the K\"ahler--Ricci flow will deform a projective variety of non-negative Kodaira dimension via a finite number of divisorial contractions and metric flips in the Gromov--Hausdorff topology, then the flow will converge to a twisted K\"ahler--Einstein metric $\omega_{\text{can}}$ on the canonical model, satisfying $$\text{Ric}(\omega_{\text{can}}) = - \omega_{\text{can}} + \omega_{\text{WP}}$$ on the regular part of the canonical model. 

It is known \cite{TianZhang2006} that the K\"ahler--Ricci flow exists on a maximal time interval $[0,T)$ given by the cohomological criterion $$T \ = \ \sup \{ t \in \mathbb{R} : [\omega_0] + t[K_M] > 0 \}.$$ This yields a sharp local existence theorem for the K\"ahler--Ricci flow. In particular, we observe that the K\"ahler--Ricci flow admits a long-time solution if and only if $K_M$ is nef.

Motivated by the abundance conjecture in birational algebraic geometry,\footnote{Which states that if $K_M$ is nef, then $K_M$ is semi-ample. The conjecture has been verified up to complex dimension 3.} we will assume that the canonical bundle $K_M$ is semi-ample, i.e., for some $\ell>0$, $K_M^{\ell}$ is basepoint free.\footnote{Explicitly, this means that for any $p \in M$, there exists $s \in H^0(M, K_M^{\ell})$ such that $s(p) \neq 0$}

There are three cases to consider concerning the long-time solutions of the K\"ahler--Ricci flow: 

If $\kappa =0$, $M$ is a Calabi--Yau manifold, and \cite{Cao} informs us that the K\"ahler--Ricci flow converges to the unique K\"ahler--Einstein metric in $[\omega_0]$. If $\kappa = m$, then $K_M$ is nef and big, and we know from \cite{TianZhang, Tsuji} that the K\"ahler--Ricci flow converges weakly to the canonical metric on $M_{\text{can}}$, which is smooth on the regular part of $M_{\text{can}}$. Further, it was shown in \cite{SongKECurrents}  that the metric completion of the K\"ahler--Einstein metric on $M_{\text{can}}^{\circ}$ is homeomorphic to $M_{\text{can}}$ itself. 

Assume now that $0 < \kappa < m$. The linear system $| K_M^{\ell} |$ (for $\ell>0$ large enough so that $K_M^{\ell}$ is basepoint free) furnishes a holomorphic map with connected fibers\footnote{For details, see, e.g., \cite{Lazarsfeld}.} $$f : M \longrightarrow N \ \hookrightarrow \ \mathbb{P}^{n_{\ell}}, \hspace*{1cm} n_{\ell} = \dim H^0(M, K_M^{\ell}) - 1,$$ whose image is a normal projective variety\footnote{Again, for $\ell>0$ sufficiently large, see, e.g., \cite{Lazarsfeld}. The fact that the image of $f$ is an analytic subvariety of $\mathbb{P}^{n_{\ell}}$ is granted to us by Remmert's proper mapping theorem \cite[p. 162]{GunningRossi}.} $N$ of (complex) dimension\footnote{See, e.g., \cite[Theorem 2.1.33]{Lazarsfeld}.} $\kappa$ which we call the canonical model of $M$.\footnote{In more detail, if $\{ s_0, ..., s_{n_{\ell}} \}$ is a basis for $H^0(M, K_M^{\ell})$, then $f$ is defined by sending $p \in M$ to the point $[s_0(p) : \cdots : s_{n_{\ell}}(p)]$ in projective space. Since $K_M^{\ell}$ is basepoint free, this map is well-defined.} Moreover, the generic fibers of $f$ are Calabi--Yau manifolds of dimension $m - \kappa$. This fibration is often referred to as an Iitaka fibration or Calabi--Yau fibration. In the intermediate Kodaira dimension case, we have: 

\vspace*{-0.1cm}

\subsection*{Conjecture 1.3.1}\label{KRFMainConjecture}
The analog of \nameref{MainConjecture} holds for the K\"ahler--Ricci flow, assuming a uniform Ricci bound on $M \backslash S$. \\

It is expected that the Ricci curvature along the flow remains uniformly bounded on $M \backslash S$ (see, e.g., \cite{TianSurvey, TianZhang}). This is known to be the case when the generic fibers are tori or finite \'etale quotients of tori \cite{GTZ2013,HeinTosattiRemarks,TosattiZhangInfiniteTimeSingularities, FongZhangKRF}. Moreover, using the higher-order estimates obtained in \cite{HeinTosatti}, the Ricci curvature was shown recently in \cite{FongLee} to be uniformly bounded when the generic fibers are all biholomorphic. In \cite{SongTianZhang2019}, part (ii) of \nameref{KRFMainConjecture} was confirmed in general, and part (iii) when $N$ has at worst orbifold singularities. Moreover, if $N$ is smooth and the codimension-one part of the discriminant locus of the Iitaka map has simple normal crossings, \nameref{KRFMainConjecture} was verified in \cite{GTZ2019}.

\subsection*{Acknowledgements}
The author is indebted to his supervisors: Ben Andrews and Gang Tian. He would also like to express his gratitude to Zhou Zhang, Behrouz Taji, Wangjian Jian, and Bin Guo for many useful discussions, as well as Valentino Tosatti for pointing out the error in the previous version of this paper. The author would like to thank the referee for pointing out the gap in the proof of Lemma 2.5, as well as bringing reference \cite{BingWang} to his attention.

\section{Proof of the Main Theorem}
Let us show how the partial second-order estimate yields the global geometric consequences, using the arguments in \cite{TosattiZhang}. To simplify the notation in this section, we write $\omega$ for the canonical metric $\omega_{\text{can}}$.

\subsection*{Theorem 2.1}\label{Hausdorff}
Let $(X, d_X)$ denote the metric completion of $(N_0, \omega)$, and set $S_X = X \backslash N_0$. Then $(X, d_X)$ is a compact length metric space and $S_X$ has real Hausdorff codimension at least 2.

 \begin{proof}
Let $\omega_{\text{cone}}$ be the cone metric with cone angles $2\pi \alpha_i$ along each component $E_i$. Since each $\alpha_i \in (0,1] \cap \mathbb{Q}$, there is an integer $m_i$ such that $1-\alpha_i < 1 - \frac{1}{m_i} < 1$. We may therefore control $\omega_{\text{cone}}$ by an orbifold K\"ahler metric $\omega_{\text{orb}}$ with orbifold order $\frac{1}{m_i}$ along $E_i$. We show that \nameref{Hausdorff} can be established from the weaker estimate: \begin{eqnarray}\label{ORBCONTROL}
\pi^{\ast} \omega & \leq & \frac{C}{| s_F |^{2\varepsilon}} \left( 1 - \sum_{i=1}^{\mu} \log | s_i |_{h_i} \right)^d \omega_{\text{orb}}.
\end{eqnarray}
Let $\text{dist}_{\omega}$ be the distance function associated to $\omega$. To show that $(X, d_X)$ is a compact length metric space, we need only show that \begin{eqnarray}\label{DIAM}
\sup_{y_1, y_2 \in N_0} \text{dist}_{\omega}(y_1, y_2)  & \leq & C.
\end{eqnarray}
The fact that $(X, d_X)$ is a length space follows from general theory. To this end, we write $$E = \bigcup_{p=1}^n E_p'$$ as a disjoint union of connected (possibly empty) $(n-p)$--dimensional relatively compact complex submanifolds of $\widetilde{N}$. These submanifolds are defined recursively: Set $E_{n+1}' = \emptyset$ and declare for multi-indices $J = (j_1, ..., j_p)$ with $1 \leq j_1, ..., j_p \leq \mu$, that $$E_p' = \bigcup_{| J | = p} (E_{j_1} \cap \cdots \cap E_{j_p})\backslash E_{p+1}'.$$ For any open neighborhood $U$ of $E_{p+1}'$, the set $E_p' \backslash U$ is compact. Now granted $\rho, \beta >0$ small, we cover $E$ by $N(\rho)$ open sets $\{ V_i(\rho) \} \subset \widetilde{N}$ such that \begin{eqnarray*}
\rho^{2n-2 + \beta} N(\rho) \longrightarrow 0, \hspace*{1cm} \text{as} \ \rho \to 0.
\end{eqnarray*}
Each $V_i(\rho)$ is contained in a chart which we identify with a unit polydisk $\Delta^n$ in $\mathbb{C}^n$. Let $(w_1, ..., w_n)$ be the coordinates in this polydisk. These coordinates can be chosen such that $E = \{ w_1 \cdots w_k =0 \}$ for some $1 \leq k \leq n$, and \begin{eqnarray*}
V_i(\rho) &=& \{ w \in \Delta^n \ : \ | w_j | < \rho^{m_j}, \ \ 1 \leq j \leq k, \ \ | w_j | < \rho, \ \ k+1 \leq j \leq n \},
\end{eqnarray*}
where for simplicity, we write $m_j$ for the orbifold order of $\omega_{\text{orb}}$ along $\{ w_j = 0 \}$. Let $q : \Delta^n \to \Delta^n$ be the map $$q(w_1, ..., w_n) = (w_1^{m_1}, ..., w_k^{m_k}, w_{k+1}, ..., w_n)$$ so that $q^{-1}(V_i(\rho))$ is the polydisk of polyradius $(\rho, ..., \rho) \in \mathbb{R}^n$, which we simply denote by $\Delta^n(\rho)$. Let $\omega_{\mathbb{C}^n}$ be the Euclidean metric on $\mathbb{C}^n$ and let $L \subset \{ 1, ..., k \}$ be the indexing set for the $\pi$--exceptional components of $E$. On $\Delta^n \backslash E$, the estimate \eqref{ORBCONTROL} reads \begin{eqnarray}\label{BASEORB}
q^{\ast} \pi^{\ast} \omega & \leq & \frac{C}{\prod_{\ell \in L} | w_{\ell} |^{2\varepsilon}} \left( 1 - \sum_{i=1}^k \log | w_i | \right)^d \omega_{\mathbb{C}^n}.
\end{eqnarray}
\subsection*{Claim:} Let $\Delta^{\ast, k}(\rho) = \{ (w_1, ..., w_k) \in \mathbb{C}^k : 0 < | w_j | < \rho, \ \ 1 \leq j \leq k \}$ denote the punctured polydisk in $\mathbb{C}^k$ of polyradius $(\rho, ..., \rho) \in \mathbb{R}^k$. For any two points $q_1, q_2 \in \Delta^{\ast,k}(\rho) \times \Delta^{n-k}(\rho) \subset \Delta^n(\rho)$, there is a path $\widetilde{\gamma} \subset \Delta^{\ast,k}(\rho) \times \Delta^{n-k}(\rho)$ connecting $q_1, q_2$ whose length is controlled in the following way: \begin{eqnarray*}
\text{length}_{q^{\ast} \pi^{\ast} \omega}(\widetilde{\gamma}) & \leq & C \rho^{1-\varepsilon} (- \log \rho)^d, 
\end{eqnarray*}
where $0 < \varepsilon < 1$. \\

The diameter estimate \eqref{DIAM} follows readily from the claim. Indeed, let $\gamma$ be the image $q(\widetilde{\gamma})$, whose $\pi^{\ast} \omega$--length coincides with the $q^{\ast}\pi^{\ast} \omega$--length of $\widetilde{\gamma}$. From the claim, any two points in $V_i(\rho) \backslash E$ can be joined by a curve $\gamma$ in $V_i(\rho)\backslash E$ with $$\text{length}_{\pi^{\ast} \omega}(\gamma) \ \leq \ C \rho^{1-\varepsilon}(-\log \rho)^d.$$ Since the open sets $\{ V_i(\rho)\}$ cover $E$, and $\rho^{1-\varepsilon}(-\log \rho)^d \to 0$ as $\rho \to 0$, this proves \eqref{DIAM}. \\

\textit{Proof of Claim:} Convert to polar coordinates $w_i = r_i \exp(\sqrt{-1} \vartheta_i)$ and translate \eqref{BASEORB} accordingly: \begin{eqnarray}\label{POLAR}
q^{\ast} \pi^{\ast} \omega & \leq & \frac{C}{\prod_{\ell \in L} r_{\ell}^{2\varepsilon}} \left( 1 - \sum_{i=1}^k \log r_i \right)^d \left( \sum_{j=1}^n dr_j^2 + r_j^2 d\vartheta_j^2 \right). 
\end{eqnarray}
Set $q_1 = (r_1 \exp(\sqrt{-1} \vartheta_1), ..., r_n \exp(\sqrt{-1} \vartheta_n))$ with $r_1, ..., r_k >0$ and $r_{k+1}, ..., r_n \geq 0$. If $r_j =0$, we set $\vartheta_j=0$. Let $\gamma_1$ be the path parametrized by $s \in [0,1]$ given by the formula $\gamma_1(s)=(\gamma_1^{(1)}(s), ..., \gamma_1^{(n)}(s))$, where $$\gamma_1^{(j)}(s) \ : = \ \left( \frac{1}{2} s \rho + (1-s) r_j \right)\exp(\sqrt{-1} \vartheta_j).$$ This gives a path in $\Delta^{\ast, k}(\rho) \times \Delta^{n-k}(\rho)$ starting at $q_1$ with endpoint on the distinguished boundary $$S(\rho/2) \ = \ \left \{ w \in \Delta^n(\rho) \ : \ | w_j | = \frac{1}{2} \rho, \ 1 \leq j \leq n \right \}.$$ The Euclidean norm of $\gamma_1'$ is at most $\rho$. We deduce from \eqref{POLAR} that \begin{eqnarray*}
\text{length}_{q^{\ast} \pi^{\ast} \omega}(\gamma_1) & \leq & C \rho \int_0^1 \prod_{\ell} ( s \rho /2 + (1-s) r_{\ell} )^{-\varepsilon} (1 - k \log(\rho s /2) )^{\frac{d}{2}} ds \\
& \leq & C \rho^{1-\varepsilon} \int_0^1 s^{-\varepsilon} (1-k\log(\rho s/2) )^{\frac{d}{2}} ds. 
\end{eqnarray*}
By Young's inequality, we have  \begin{eqnarray*}
\int_0^1 s^{-\varepsilon} (1-k\log(\rho s/2))^{\frac{d}{2}} ds & \leq & \int_0^1 s^{-2\varepsilon} ds + \int_0^1 (1-k\log(\rho s/2))^d ds \\
& \leq & \frac{1}{1-2\varepsilon} + C(-\log \rho)^d.
\end{eqnarray*}
Hence, we see that \begin{eqnarray}\label{DESIRED}
\text{length}_{q^{\ast} \pi^{\ast} \omega}(\gamma_1) & \leq & C \rho^{1-\varepsilon} (- \log \rho)^d.
\end{eqnarray}
Similarly, with the distinguished boundary $S(\rho/2)$ diffeomorphic to the real $n$--torus $\mathbb{T}^n = \mathbb{R}^n/\mathbb{Z}^n$, the distinguished boundary $S(\rho/2)$ is connected and \eqref{DESIRED} implies \begin{eqnarray*}
\text{diam}_{q^{\ast} \pi^{\ast} \omega}(S(\rho/2)) & \leq & C \rho^{1-\varepsilon} (-\log \rho)^d. 
\end{eqnarray*}
This proves the claim and therefore the first part of the theorem.\\

We now show that the real Hausdorff codimension of $S_X$ is at least 2. Recall from \cite{TosattiWeinkoveYang, GTZ2013} that there is a local isometric embedding $\Phi : (N_0, \omega) \longrightarrow (X, d_X)$ with open dense image $X_0 \subset X$. Density is obtained by proving that $\nu(S_X) =0$, where $\nu$ is the renormalized limit measure (see \cite{CC1}). For any $\rho > 0$, the density of $X_0$ gives us a cover of $S_X$: \begin{eqnarray}\label{COVER}
S_X \subset \bigcup_{i=1}^{N(\rho)} \overline{U_i(\rho)}, \hspace*{1cm} U_i(\rho) = \Phi(\pi(V_i(\rho)) \cap N_0).
\end{eqnarray}
The length space structure of $(X, d_X)$ tells us that for each $i$, \begin{eqnarray*}
\text{diam}_{d_X} \overline{U_i(\rho)} &=& \text{diam}_{d_X} U_i(\rho) \ = \ \sup_{p,q \in U_i(\rho)} \inf_{\eta} \text{length}_{d_X}(\Gamma), 
\end{eqnarray*}
where the infimum is taken over all curves $\Gamma$ in $X$ which connect $p$ and $q$. From the proof of the first part of the theorem, however, $p,q$ can be joined by a curve of the form $\Phi(\pi(\gamma))$, where $\gamma$ lies in $V_i(\rho)\backslash E$ and whose length is controlled in the desirable way \eqref{DESIRED}. With $\Phi$ a local isometry, it follows that \begin{eqnarray*}
\text{length}_{d_X}(\Phi(\pi(\gamma)) \ = \ \text{length}_{\omega}(\pi(\gamma)) \ = \ \text{length}_{\pi^{\ast} \omega} (\gamma),
\end{eqnarray*}
and so \begin{eqnarray}\label{DIAM}
\text{diam}_{d_X} \overline{U_i(\rho)} \ \leq \ C \rho^{1-\varepsilon} (-\log \rho)^d.
\end{eqnarray} Fix $\beta, \delta >0$ small, and for any $\eta >0$, choose $\rho > 0$ sufficiently small such that $$C \rho^{1-\varepsilon} (-\log \rho )^d < \eta.$$ Let $v_{2n}$ denote the volume of the unit ball in $\mathbb{R}^{2n}$. Then \begin{eqnarray*}
\mathcal{H}_{d_X, \eta}^{2n - 2 + \beta + \delta}(S_X) & \leq & \sum_{i=1}^{N(\rho)} v_{2n} \text{diam}_{d_X}^{2n - 2 + \beta + \delta} ( \overline{\psi(\pi(V_i(\rho)) \cap N_0)}) \\
& \leq & C N(\rho) \left[ \rho^{1-\varepsilon} ( -\log \rho)^d \right]^{2n - 2 + \beta + \delta} \\
& \leq & C \left( N(\rho) \rho^{2n - 2 + \beta} \right)  \rho^{\delta - \varepsilon(2n-2+\beta+\delta)} (-\log \rho)^{d(2n - 2 + \beta + \delta)} \\
& \to & 0,
\end{eqnarray*}
as $\rho \to 0$. Note that here we choose $\varepsilon >0$ such that $\delta - \varepsilon(2n-2+\beta + \delta) >0$. Therefore, \begin{eqnarray*}
\mathcal{H}_{d_X}^{2n - 2 + \beta + \delta}(S_X) &=& \lim_{\eta \to 0} \mathcal{H}_{d_X, \eta}^{2n - 2 + \beta + \delta}(S_X) \ = \ 0,
\end{eqnarray*}
for any small $\beta, \delta >0$. This shows that $\dim_{\mathcal{H}}(S_X) \leq 2n -2$ and completes the proof.
\end{proof}

\subsection*{Corollary 2.2}
As $t \to 0$, the metric spaces $(M, \text{dist}_{\widetilde{\omega}_t})$ converge in the Gromov--Hausdorff topology to $(X,d_X)$, where $(X,d_X)$ is the metric completion of $(N_0, \omega)$. \begin{proof}
The argument in \cite[p. 8--9]{TosattiZhang} can be followed verbatim.
\end{proof}

\subsection*{Theorem 2.3}\label{GHLimit}
The Gromov--Hausdorff limit $(X, d_X)$ is homeomorphic to $(N, \omega)$. \\

Again, we will show how \nameref{MainTheorem} implies \nameref{GHLimit}. From \cite{TosattiWeinkoveYang,GTZ2013}, we know that the Gromov--Hausdorff limit $(X, d_X)$ is isometric to the metric completion of $(N_0, \omega)$. We let $\Phi : (N_0, \omega) \longrightarrow (X, d_X)$ denote the isometric embedding. The Yau Schwarz lemma proved in \cite{TosattiAdiabaticLimits} informs us that for any K\"ahler metric $\omega_N$ on $N$, we have $\widetilde{\omega}_t \geq C^{-1} f^{\ast} \omega_N$. In particular, this yields a uniform (independent of $t$) bound on the Lipschitz constant for $f : (M, \widetilde{\omega}_t) \longrightarrow (N, \omega_N)$. Passing to a subsequence $t_i \searrow 0$, we get a surjective Lipschitz map $$h : (X, d_X) \longrightarrow (N, \omega_N),\hspace*{1cm} h \circ \Phi = \text{id}.$$ To indicate the proof of \nameref{GHLimit}, let $\pi : \widetilde{N} \longrightarrow N$ be the composition of blow-ups with smooth centers such that $\widetilde{N}$ is smooth and $E = \pi^{-1}(D)$ is an snc divisor. Assume there is a continuous surjective map $p : \widetilde{N} \longrightarrow X$ such that $\pi = h \circ p$. If $\pi = \text{id}$ (the case considered in \cite{GTZ2019}), then $h \circ p = \text{id}$ and $p$ is necessarily injective, as required.\footnote{In more detail, $p$ is a continuous bijection between compact Hausdorff spaces, which implies that $p$ is a homeomorphism.} In general, however, $\pi$ fails to be the identity map, and as observed in \cite{TosattiZhang}, it suffices to prove that $p$ factors through $\pi$. Granted this, there exists a continuous surjective map $\overline{p} : N \longrightarrow X$ such that $p = \overline{p} \circ \pi$. Hence, $\pi = h \circ \overline{p} \circ \pi$ and since $\pi$ is surjective, $h \circ \overline{p} = \text{id}$ and we see that $\overline{p}$ is the desired homeomorphism. \nameref{GHLimit} therefore follows from the following two lemmas:

\subsection*{Lemma 2.4}\label{SURJ}
There is a continuous surjective map $p : \widetilde{N} \longrightarrow X$ such that $\pi = h \circ p$. \begin{proof}
The map $p : = \Phi \circ \pi : \pi^{-1}(N_0) \longrightarrow \Phi(N_0) \subset X$ is a homeomorphism. So we want to extend $p$ to a map $p : \widetilde{N} \longrightarrow X$. To define the extension, let $y \in \widetilde{N} \backslash \pi^{-1}(N_0)$. We choose a sequence $y_j \in \pi^{-1}(N_0)$ such that $y_j \to y$ in $\widetilde{N}$. For $j$ sufficiently large, the $y_j$ are all contained in some sufficiently small polydisk $\Delta^n(\rho)$, so for all $m \geq 0$, \begin{eqnarray*}
\text{dist}_{\pi^{\ast} \omega} (y_j, y_{m+j}) & \leq & C \rho^{1-\varepsilon} (-\log \rho)^d \ \to \ 0.
\end{eqnarray*}
The sequence $\pi(y_j)$ is therefore Cauchy in $(N_0, \omega)$ and converges to a unique limit point $\overline{y}$ in $X$. We then simply define $p$ to send $y$ to $\overline{y}$. It can be readily verified that this definition is well-defined. To see that $p$ is surjective, take $\overline{y} \in X$, and let $ y_j  \in (N_0, \omega)$ be the Cauchy sequence converging to this point. Choose any preimages $y_j' = \pi^{-1}(y_j) \in \widetilde{N}$. This sequence will admit a convergent subsequence by the compactness of $\widetilde{N}$, limiting to some $y' \in \widetilde{N}$. We then necessarily have that $p(y') = \overline{y}$. Finally, if $y' = \pi(y) \in N$, it follows that $\pi(y_j) = h(p(y_j)) \to y'$ and so $h(\overline{y}) = y'$, which completes the proof. 
\end{proof}

\subsection*{Lemma 2.5}\label{FACTOR}
For $x,y \in \widetilde{N}$ with $\pi(x) = \pi(y)$, we have $p(x) = p(y)$. 

\begin{proof}
Assume $\pi(x) = \pi(y)$ for two points $x,y \in E$ and let $x_k, y_k \in \pi^{-1}(N_0)$ be sequences such that $x_k \to x$ and $y_k \to y$. For any small neighbourhood $\Delta^n(x, \rho)$ of $x$ and any small neighbourhood $\Delta^n(y,\rho)$ of $y$, we may choose $k$ sufficiently large such that $x_k \in \Delta^n(x, \rho)$ and $y_k \in \Delta^n(y, \rho)$. Given any $\ell \geq 0$, the partial second-order estimate affords the estimates: $$\text{dist}_{\pi^{\ast} \omega}(x_k, x_{k+\ell}) \leq C \rho^{1-\varepsilon}(-\log \rho)^d. \hspace{1cm}  \text{dist}_{\pi^{\ast} \omega}(y_k, y_{k+\ell}) \leq C \rho^{1-\varepsilon}(-\log \rho)^d.$$   

The sequences $(x_k), (y_k)$ are Cauchy with respect to $\text{dist}_{\pi^{\ast}\omega}$ and therefore converge to points $p(x), p(y) \in S_X$ in the singular part of the metric completion $(X, d_X)$. Recall that $S_X$ admits a covering of the form: $$S_X = \bigcup_{i=1}^{N(\rho)} \overline{U_i(\rho)}, \hspace{1cm} U_i(\rho) = \Phi(V_i(\rho) \cap N_0),$$ where $\text{diam}_{d_X}(\overline{U_i(\rho)}) \leq C \rho^{1- \varepsilon} (-\log \rho)^d$. Assume $p(x) \in \overline{U_i(\rho)}$ and $p(y) \in \overline{U_j(\rho)}$ for some fixed indices $i,j$. From the diameter control of $S_X$, and the density of the regular part of $X$, choose $\widetilde{x}, \widetilde{y} \in X \backslash S_X$ such that $p(x)$ and $\widetilde{x}$ can be connected by a curve $\gamma_x$ with $\text{length}_{d_X}(\gamma_x) \leq C \rho^{1-\varepsilon}(-\log \rho)^d$, and $p(y)$ connects to $\widetilde{y}$ via a curve $\gamma_y$ with $\text{length}_{d_X}(\gamma_y) \leq C \rho^{1-\varepsilon}(-\log \rho)^d$.  We showed previously that any two points in $V_i(\rho)\backslash E$ may be joined by a curve $\gamma$ which does not intersect $E$ and whose length is bounded from above by a constant multiple of $\rho^{1-\varepsilon}(-\log \rho)^d$.  Since $\widetilde{x}, \widetilde{y} \in X \backslash S_X \simeq N_0 \simeq \widetilde{N} \backslash E$, write $\widetilde{x} = \Phi(\pi(\hat{x}))$, $\widetilde{y} = \Phi(\pi(\hat{y}))$ for two points $\hat{x}, \hat{y}$ in $V_k(\rho) \backslash E$, for some $k$, and therefore estimate: \begin{eqnarray*}
d_X(p(x),p(y)) & \leq & d_X(p(x), \widetilde{x}) + d_X(\widetilde{x}, \widetilde{y}) + d_X(\widetilde{y}, p(y)) \\
& \leq & C \rho^{1-\varepsilon}(-\log \rho)^d + d_X(\Phi(\pi(\hat{x}), \Phi(\pi(\hat{y}))) \\
&=& C \rho^{1-\varepsilon}(-\log \rho)^d + \text{dist}_{\pi^{\ast} \omega}(\hat{x}, \hat{y}) \\
& \leq & C' \rho^{1-\varepsilon}(-\log \rho)^d. 
\end{eqnarray*}

Hence, $\pi(x) = \pi(y)$ implies $p(x)=p(y)$, as required.
\end{proof}

\subsection*{Remark 2.6}
Note that the proof of \nameref{FACTOR} in \cite{TosattiZhang} requires the machinery of special K\"ahler geometry and does not follow from the (partial) second-order estimate \eqref{PARTIALSECOND}. The proof presented here, however, explicitly uses the partial second-order estimate and uses less specific machinery.

\hfill

\scshape{Mathematical Sciences Institute, Australian National University, Acton, ACT 2601, Australia}

\scshape{BICMR, Peking University, Beijing, 100871, People's republic of china}

\textit{E-mail address}: \texttt{kyle.broder{@}anu.edu.au}

\end{document}